\documentclass[conference]{IEEEtran}
\IEEEoverridecommandlockouts
\usepackage{cite}
\usepackage{amsmath,amssymb,amsfonts}
\usepackage{makecell, multirow, threeparttable, tabularx, array, arydshln}
 
\usepackage{hyperref} 
\usepackage{booktabs} 
\usepackage{cleveref}
\crefname{equation}{Eq.}{Eq.}
\crefname{figure}{Fig.}{Fig.}
\crefname{table}{Table}{Table}
\crefname{algorithm}{Alg.}{Algorithms}
\Crefname{ALC@unique}{Line}{Lines}

\usepackage{amsthm}
\usepackage{algorithm}
\usepackage{algorithmic}
\usepackage{graphicx}

\newtheorem{thm}{Theorem}
\newtheorem{lem}{Lemma}

\usepackage{textcomp}
\usepackage{xcolor}
\usepackage[a4paper, total={184mm,239mm}]{geometry}
\def\BibTeX{{\rm B\kern-.05em{\sc i\kern-.025em b}\kern-.08em
    T\kern-.1667em\lower.7ex\hbox{E}\kern-.125emX}}
    
\begin{document}

\title{Computing Effective Resistances on Large Graphs Based on Approximate Inverse of Cholesky Factor
\thanks{This work is supported by NSFC under grant No. 62090025, and Beijing Science and Technology Plan (No. Z221100007722025).}
}

	\author{\IEEEauthorblockN{Zhiqiang Liu, Wenjian Yu}
		\IEEEauthorblockA{Dept. Computer Science \& Tech., BNRist, Tsinghua University, Beijing 100084, China\\
			Email: liu-zq20@mails.tsinghua.edu.cn, yu-wj@tsinghua.edu.cn
	}}

\maketitle

\begin{abstract}
Effective resistance, which originates from the field of circuits analysis, is an important graph distance in spectral graph theory. It has found numerous applications in various areas, such as graph data mining, spectral graph sparsification, circuits simulation, etc. However, computing effective resistances accurately can be intractable and we still lack efficient methods for estimating effective resistances on large graphs. In this work, we propose an efficient algorithm to compute effective resistances on general weighted graphs, based on a sparse approximate inverse technique. Compared with a recent competitor, the proposed algorithm shows several hundreds of speedups and also one to two orders of magnitude improvement in the accuracy of results. Incorporating the proposed algorithm with the graph sparsification based power grid (PG) reduction framework, we develop a fast PG reduction method, which achieves an average 6.4X speedup in the reduction time without loss of reduction accuracy. In the applications of power grid transient analysis and DC incremental analysis, the proposed method enables 1.7X and 2.5X speedup of overall time compared to using the PG reduction based on accurate effective resistances, without increase in the error of solution. 
\end{abstract}

\begin{IEEEkeywords}
Effective resistances, power grid reduction, DC incremental analysis, transient analysis.
\end{IEEEkeywords}

\section{Introduction}

Effective resistance is an important metric that measures the vertex similarity in a graph. It has found tremendous applications in a variety of areas, including graph data mining \cite{SCWWW15,SCIJCAI16,SCKDD21}, spectral graph sparsification \cite{EffRes08,kmp,fegrass,liu2022pursuing} and circuit simulation \cite{ResMORICCAD14,ResMORDAC19,EffResTCAS20}, etc. Computing effective resistances for many node pairs on large graphs is a computationally challenging task. There are several methods for estimating effective resistances in the literature. A random projection based method was introduced in \cite{SCWWW15}, but it still takes a huge amount of time to compute effective resistances with high accuracy. The methods based on random walk or random spanning tree generation were proposed in \cite{SCIJCAI16,SCKDD21}, but they can only handle unweighted graphs. A method based on infinity mirror techniques was presented in \cite{EffResTCAS20}, but it is under the assumption that the graph is of two-dimensional grid structure. Despite of the  importance of effective resistances, we still lack efficient methods for computing effective resistances on general large graphs.

The goal of power grid (PG) reduction is to reduce the original large power grid to a smaller one which can preserve the electrical behavior of port nodes. There are mainly three types of PG reduction methods in the literature: moment matching based methods \cite{mmmor1,mmmor2}, node elimination based methods \cite{ticer,SIPICCAD08} and multigrid-like methods \cite{amgmor1,amgmor2}. Moment matching based methods \cite{mmmor1,mmmor2} cannot reduce power grids efficiently because there are hundreds of thousands of port nodes in typical power grids. Node elimination based methods \cite{ticer,SIPICCAD08} have been successful for reducing tree-like RC networks, but they tend to generate much denser models with even more edges than original models when dealing with mesh-like power grids. The multigrid-like methods \cite{amgmor1,amgmor2} can generate realizable and sparse reduced models, but their error is hard to control.

Effective resistances have been utilized for developing more efficient PG reduction methods \cite{ResMORICCAD14,ResMORDAC19}. They employ the effective resistances based sampling approach \cite{EffRes08} to sparsify the dense reduced models. The method proposed in \cite{ResMORDAC19} does not scale well to large problems due to the high computational complexity, as admitted in \cite{ResMORDAC19}. The method proposed in \cite{ResMORICCAD14} is more scalable since it leverages the techniques of graph partitioning and effective resistances based port merging. However, it computes effective resistances accurately, which still takes a large amount of time for large-scale power grids.

In this paper, we aim to develop an efficient algorithm for computing effective resistances on large graphs and also a fast PG reduction method. Our main contributions are summarized as follows.

1) We propose an efficient algorithm for computing effective resistances on general weighted graphs, which is based on a sparse approximate inverse technique for the Cholesky factor of Laplacian matrix.

2) Incorporating the proposed algorithm for computing effective resistances with the PG reduction framework proposed in \cite{ResMORICCAD14}, we develop a fast PG reduction method.

Extensive experiments have been conducted to validate the efficiency and the accuracy of the proposed algorithm for computing effective resistances, which shows an average 168X speedup and also significant reduction in errors over the recent competitor \cite{SCWWW15}. The proposed algorithm is highly scalable. For a graph with 6.0E7 nodes and 1.0E8 edges, effective resistances of all edges can be computed within about 8 minutes, with an estimated average relative error of only 0.17\%. It also derives a fast PG reduction method, which achieves an average 6.4X speedups over the method based on accurate effective resistances, with no increase in reduction errors. The resulted fast PG reduction method can be utilized in many downstream applications. For example, it brings 1.7X and 2.5X speedups of overall time for power grid transient analysis and DC incremental analysis, respectively.

\section{Background}
\subsection{Effective Resistances and Their Applications}
	Suppose $G =(V, E, w)$ denotes a weighted undirected graph, where $V$ and $E$ are the sets of vertices (nodes) and edges, $w$ is a positive weight function. Let $n=|V|$ and $m=|E|$. 
The incidence matrix $B$ is defined to be an $m \times n$ matrix such that each row of $B$ corresponds to an edge in $E$ and each column of $B$ corresponds to a node in $V$. The entries of $B$ satisfy:
	\begin{equation}
	B(e,v)=\left\{
	\begin{aligned}
	&1, \quad v ~ is ~ e's ~ head \\
	&-1, \quad v ~ is ~ e's ~ tail \\
	&0, \quad \textrm{otherwise} ~ .
	\end{aligned}
	\right.
	\end{equation}

Let $W$ denote an $m \times m$ diagnoal matrix with $W(e,e)=w(e)$. Then the Laplacian matrix $L_G \in \mathbb{R}^{n \times n}$ is defined as:
\begin{equation}
L_G=B^TWB~.
\end{equation}

Given two nodes $p$ and $q$, the \textit{effective resistance} $R_G(p,q)$ across $p$ and $q$ refers to the voltage difference between $p$ and $q$ when a unit current flows into $G$ through $p$ and leaves through $q$. Formally, it can be defined as:
\begin{equation}
\label{equ:effresdef}
R_G(p,q)=e_{p,q}^TL_G^{\dagger}e_{p,q}~.
\end{equation}
Here $L_G^{\dagger}$ denotes pseudo-inverse of $L_G$ and $e_{p,q}=e_p-e_q$, where $e_p$ is the $p$-th column of the identy matrix. $L_G$ is singular as the smallest eigenvalue is $0$. To handle the singularity of $L_G$, we introduce a ground node and connect it to a randomly selected node in each connected component of $G$. It corresponds to adding some small positive values to the diagonal elements of $L_G$. To simplify the notations, we still use $L_G$ to denote the resulted symmetric diagonally dominant (SDD) matrix. 

Below we briefly review the power grid reduction method proposed in \cite{ResMORICCAD14}, which employs effective resistances based graph sparsification.


  
Power grid reduction \cite{EffResMixTime} aims to reduce the original large power grid to a small one which contains all the port nodes. A port node is defined to be a node which is connected to a voltage or current source. The voltage drops of port nodes can be obtained by analyzing the reduced model, enabling more efficient analysis for large power grids. Schur complement based method \cite{SIPICCAD08} is adopted to eliminate the non-port nodes without loss of accuracy but tends to generate much denser models. To address this issue, a graph sparsification based power grid reduction approach was introduced \cite{ResMORICCAD14}, which consists of circuits partitioning, Schur complement based reduction, effective resistances based port merging and effective resistances based graph sparsification. We remark that all the port nodes have important physical information and should be preserved. Although at least half of the port nodes were eliminated in \cite{ResMORICCAD14}, the algorithm can be easily extended to the case where all the port nodes should be kept. We summarize the modified version as \cref{alg:ResMOR}.  
	\begin{algorithm}[h]
		\caption{Power Grid Reduction via Effective Resistances Based Graph Sparsification \cite{ResMORICCAD14}}
		\label{alg:ResMOR}
		\begin{algorithmic}[1]
			\setcounter{ALC@unique}{0}
			\REQUIRE The original power grid. 
			\ENSURE The reduced power grid.
			\STATE Partition the original power grid into many blocks. The nodes are classified into three types: port nodes, non-port interface nodes and non-port interior nodes.
			\STATE For each block, eliminate the non-port interior nodes using Schur complement based method.
			\STATE For each reduced block, compute effective resistances exactly using 
   (\ref{equ:effresdef}) for all the edges.
			\STATE For each reduced block, merge the nodes based on effective resistances and then sparsify the reduced grid using effective resistances based sampling approach.
			\STATE Stitch all the reduced and sparsified models together to form the final reduced power grid.
		\end{algorithmic}
	\end{algorithm} 

Among all these steps, computing effective resistances is the most time-consuming one. So, to obtain a fast power grid reduction algorithm, more efficient methods for computing effective resistances are demanded.  

\subsection{Methods for Computing Effective Resistances}

Computing effective resistances accurately is computational challenging. For each query $(p,q)$, computing $R_G(p,q)$ requires solving linear equations whose coefficient matrix is $L_G$. It should be noted that many applications require computing effective resistances for every edge $(p,q) \in E$. Solving $|E|$ linear equations takes at least $\Omega(|E|^2)$ time, which can be prohibitive for large-scale problems.

To compute effective resistances more efficiently, a random projection based method was proposed in \cite{EffRes08}. It is based on the fact that the effective resistance $R_G(p,q)$ can be written as the distance between two vectors:
\begin{equation}
\begin{aligned}
R_G(p,q) & =e_{p,q}^TL_G^{\dagger}e_{p,q} = e_{p,q}^TL_G^{\dagger}L_GL_G^{\dagger}e_{p,q} \\
&= e_{p,q}^TL_G^{\dagger}B^TWBL_G^{\dagger}e_{p,q} \\
&= \|W^{1/2}BL_G^{\dagger}e_p-W^{1/2}BL_G^{\dagger}e_q\|_2^2 ~.\\
\end{aligned}
\end{equation}
For each node $p$, $W^{1/2}BL_G^{\dagger}e_p$ is an $m$-dimensional vector. With the Johnson-Lindenstraus Lemma, these vectors can be projected into a $k$-dimensional space, such that:
\begin{equation}
\label{equ:JLRes}
R_G(p,q) \approx \|QW^{1/2}BL_G^{\dagger}e_p-QW^{1/2}BL_G^{\dagger}e_q\|_2^2~, 
\end{equation}
where $k=O(logm)$ and $Q \in R^{k \times m}$ is a random matrix whose elements are uniformly sampled from $\pm \frac{1}{\sqrt{k}}$. 
Note that using (\ref{equ:JLRes}), only $k$ linear equations need to be solved to construct the matrix $QW^{1/2}BL_G^{\dagger}$ and then each effective resistance query can be answered in $O(k)=O(logm)$ time. A more practical version of this algorithm was presented in \cite{SCWWW15}, which incorporates the random projection based method and fast linear equation solver \cite{CMG11}. 
However, due to the big hidden constants, it still takes a huge amount of time to compute effective resistances on large graphs with reasonably high accuracy.

\section{Efficiently Computing Effective Resistances Based on Approximate Inverse of Cholesky Factor}

\subsection{The Idea}

Suppose $L_G$ is factorized with Cholesky factorization:
\begin{equation}
L_G=LL^T~,
\end{equation}
where $L$ is a lower triangular matrix. Then the effective resistance of query $(p,q)$ can be written as:
\begin{equation}
\label{equ:effres1}
\begin{aligned}
R_G(p,q) & =e_{p,q}^TL_G^{-1}e_{p,q}=e_{p,q}^TL^{-T}L^{-1}e_{p,q}\\
&=\|L^{-1}e_{p,q}\|_2^2=\|L^{-1}e_{p}-L^{-1}e_{q}\|_2^2 ~.\\
\end{aligned}
\end{equation} 
\cref{equ:effres1} shows that the effective resistance $R_G(p,q)$ can be written as the distance between the $p$-th column and the $q$-th column of $L^{-1}$. If $L^{-1}$ is available, the effective resistances can be computed efficiently, but computing and storing the dense $L^{-1}$ explicitly can be prohibitive for large-scale problems. So, our idea is to derive a sparse and approximate inverse of $L$.

\subsection{Sparse Approximate Inverse of Cholesky Factor}
Based on the observation that most elements in $L^{-1}$ are very small, we develop a method for computing a sparse approximation of $L^{-1}$. Let $Z=L^{-1}$ and $z_j$ be the $j$-th column of $Z$. Recall that $L$ is the Cholesky factor of $L_G$, it can be shown that all the diagonal elements in $L$ are positive and all the off-diagonal elements in $L$ are nonpositive \cite{directbook}.

The matrix $Z$ has some useful structural properties which are summarized as the following lemma. Here a matrix is nonnegative means that all the elements in it are nonnegative.
\begin{lem}
Suppose $Z=L^{-1}$, where $L$ is the Cholesky factor of Laplacian matrix. Then, $Z$ is nonnegative and the columns of $Z$ satisfy:
\begin{equation}
\label{equ:prop1}
z_j=\frac{1}{L_{j,j}}e_j+\sum_{i>j\&L_{i,j} \ne 0}\frac{-L_{i,j}}{L_{j,j}}z_i~,
\end{equation}
\end{lem}

\begin{proof}
\cref{equ:prop1} can be derived by multiplying $L$ on both sides. Suppose $z_i$ is nonnegative. Since $L_{i,j} \le 0$ and $L_{j,j} > 0$, \cref{equ:prop1} implies $z_j$ is nonnegative. So it can be proved by a simple mathematical induction that $Z$ is nonnegative.
\end{proof}  

\cref{equ:prop1} suggests that we can first compute the $n$-th column of $L^{-1}$ and then compute the ($n$-1)-th, ($n$-2)-th, $\cdots$, and the 1st columns one by one. Suppose we have some sparse approximations to $z_i$, denoted by $\tilde{z}_i$. Then $z_j$ can be computed approximately by:
\begin{equation}
z_j \approx z_j^* = \frac{1}{L_{j,j}}e_j+\sum_{i>j\&L_{i,j} \ne 0}\frac{-L_{i,j}}{L_{j,j}}\tilde{z}_i~.
\end{equation} 
The $z_j^*$ can be computed efficiently because the $\tilde{z}_i$s are sparse. Note that many elements in $z_j^*$ are very small, which can be set to $0$. Here we adopt a prunning strategy to control the error in the sense of $1$-norm. Note that for a vector $x$, the $1$-norm of $x$, denoted by $\|x\|_1$, is defined as the sum of absolute value of each element. Let $trunc_k(x)$ denote the vector obtained by setting the $k$ smallest elements (in the sense of absolute values) in $x$ to $0$. Our strategy is to find the largest $k$ and to set $\tilde{z}_j=trunc_k(z_j^*)$, such that:
\begin{equation}
\frac{\|\tilde{z}_j-z_j^*\|_1}{\|z_j^*\|_1} \le \epsilon~,
\end{equation} 
where $\epsilon$ is a user-defined threshold. It can be implemented by first sorting the nonzero elements in $z_j^*$ and then finding the largest $k$ which satisfies the above constraint. We summarize the algorithm for computing sparse approximate inverse of Cholesky factor as \cref{alg:iinv}.

	\begin{algorithm}[H]
		\caption{Sparse Approximate Inverse of Cholesky Factor}
		\label{alg:iinv}
		\begin{algorithmic}[1]
			\setcounter{ALC@unique}{0}
			\REQUIRE Cholesky factor of $L_S$: $L$, a user-defined threshold $\epsilon$. 
			\ENSURE A sparse approximation to $L^{-1}$: $\tilde{Z}$.
			\FOR{$j=n$ to $1$}
			\STATE Compute $z_j^* = \frac{1}{L_{j,j}}e_j+\sum_{i >j \& L_{i,j} \ne 0}\frac{-L_{i,j}}{L_{j,j}} \tilde{z}_i ~$.
			\IF{$nnz(z_j^*) \le \log n $}
			\STATE $\tilde{z}_j=z_j^*$.
			\STATE Continue.
			\ENDIF
			\STATE Find the largest $k$ such that $\frac{\|trunc_k(z_j^*)-z_j^*\|_1}{\|z_j^*\|_1} \le \epsilon~$. Set $\tilde{z}_j=trunc_k(z_j^*)$.	
			\ENDFOR
		\end{algorithmic}
	\end{algorithm}

Now we analyze the errors caused by approximating $L^{-1}$ with $\tilde{Z}$. Formally, we give a theorem below, which relates the approximation error $\|z_p-\tilde{z}_p\|_1$ to the depth of node $p$ in the filled graph. The filled graph of $G$ refers to the undirected graph corresponding to the matrix $L$ \cite{directbook}. Let $G_L=(V,F)$ denote the filled graph, where $F=\{(i,j)|i \ne j ~and~ L_{i,j} \ne 0\}$. The depth of node $p$, denoted by $depth(p)$, is defined as:
\begin{equation}
\label{equ:depth}
	depth(p)=\left\{
	\begin{aligned}
	& 0, \quad L(p+1:n,p)=\mathbf{0} \\
	& 1+\max_{i>p\&L(i,p) \ne 0}depth(i), \quad \textrm{otherwise} ~.\\
	\end{aligned}
	\right.
	\end{equation} 

\begin{thm}
		\label{thm:mainthm}
		Suppose $L$ and $\tilde{Z}=[\tilde{z}_1,\tilde{z}_2,...,\tilde{z}_n]$ are the input and output of \cref{alg:iinv}. $Z=L^{-1}=[z_1,z_2,...,z_n]$. The $depth(p)$ is defined as (\ref{equ:depth}). Then, for any node $p$:
		\begin{equation}
		\label{equ:errbnd}
		\frac{\|z_p-\tilde{z}_p\|_1}{\|z_p\|_1} \le depth(p) \times \epsilon~. 
		\end{equation}
	\end{thm}
	\begin{proof}
		We prove the theorem by induction. First note that, for node $q$ with $L(q+1:n,q)=\mathbf{0}$ or $depth(q)=0$, we have:
		\begin{equation}
		z_q=z_q^*=\tilde{z}_q=\frac{1}{L_{q,q}}e_q~,
		\end{equation}
		so $\|z_q-\tilde{z}_q\|_1=0$, which satisfies (\ref{equ:errbnd}). Now we consider node $p$ with $L(p+1:n,p) \ne \mathbf{0}$. It is hypothesized that for node $i$ with $i>p ~ and ~ L(i,p) \ne 0$, we have $\frac{\|z_i-\tilde{z}_i\|_1}{\|z_i\|_1} \le depth(i) \times \epsilon$.

		The error $\|z_p-\tilde{z}_p\|_1$ can be divided into two parts:
		\begin{equation}
		\|z_p-\tilde{z}_p\|_1=\|(z_p-z_p^*)+(z_p^*-\tilde{z}_p)\|_1 \le \|z_p-z_p^*\|_1+\|z_p^*-\tilde{z}_p\|_1 ~.
		\end{equation}
		For the second part, \cref{alg:iinv} ensures that $\|z_p^*-\tilde{z}_p\|_1 \le \epsilon \times \|z_p^*\|_1$. Recall that $z_p$ is a vector with all elements nonnegative and $z_p^*$ is a truncated version of $z_p$, it is easy to show $\|z_p^*\|_1 \le \|z_p\|_1$, so we have:
		\begin{equation}
		\|z_p^*-\tilde{z}_p\|_1 \le \epsilon \times \|z_p\|_1~.
		\end{equation}

		Now consider the first part. Recall that:
		\begin{equation}
		z_p^* = \frac{1}{L_{p,p}}e_p+\sum_{i >p \& L_{i,p} \ne 0}\frac{-L_{i,p}}{L_{p,p}} \tilde{z}_i ~.
		\end{equation}
		Here $L_{i,p} < 0$, $L_{p,p} > 0$ and $\sum_i \frac{-L_{i,p}}{L_{p,p}} \le 1$. In the rest of this proof, we use $\sum_i$ to replace $\sum_{i >p \& L_{i,p} \ne 0}$. Then
		\begin{equation}
		\|z_p-z_p^*\|_1=\|\sum_{i}\frac{-L_{i,p}}{L_{p,p}} (z_i-\tilde{z}_i)\|_1 \le \sum_{i} \frac{-L_{i,p}}{L_{p,p}} \|(z_i-\tilde{z}_i)\|_1 .
		\end{equation}
		By the inductive hypothesis, we have
		\begin{equation}
		\label{equ:prf1}
		\begin{aligned}
		\|z_p-z_p^*\|_1 & \le \sum_{i} \frac{-L_{i,p}}{L_{p,p}} \|z_i\|_1 \times depth(i) \times \epsilon \\ 
		& \le \max_i depth(i) \times \epsilon \times \sum_{i} \frac{-L_{i,p}}{L_{p,p}} \|z_i\|_1 ~.\\
		\end{aligned}
		\end{equation}
		Note that for nonnegative vectors $z_i$s, we have 
		\begin{equation}
		\label{equ:prf2}
		\sum_{i} \frac{-L_{i,p}}{L_{p,p}} \|z_i\|_1 = \|\sum_{i} \frac{-L_{i,p}}{L_{p,p}} z_i\|_1 \le \|z_p\|_1~.
		\end{equation}	
		Substituting (\ref{equ:prf2}) into (\ref{equ:prf1}), we obtain
		\begin{equation}
		\|z_p-z_p^*\|_1 \le \max_{i} depth(i) \times \epsilon \times \|z_p\|_1~.
		\end{equation}
		So
		\begin{equation}
		\frac{\|z_p-\tilde{z}_p\|_1}{\|z_p\|_1} \le (\max_{i} depth(i)+1) \times \epsilon = depth(p) \times \epsilon~. 
		\end{equation}
	
		By mathematical induction, \cref{equ:errbnd} is true for all nodes. This ends the proof. 
	\end{proof}

Note that for most graphs that stem from real world, the maximum node depth in the filled graph is not very large, as reported in the next section. By setting a sufficiently small $\epsilon$, $\tilde{Z}$ approximates $L^{-1}$ very well.

\subsection{The Overall Algorithm and Discussion}
Based on the approximate inverse technique, the effective resistance $R_G(p,q)$ can be computed as:
\begin{equation}
\label{equ:approxeffres}
R_G(p,q)=\|L^{-1}e_{p}-L^{-1}e_{q}\|_2^2 \approx \|\tilde{z}_p-\tilde{z}_q\|_2^2 ~.
\end{equation}

 Now we analyze the errors of effective resistances. Let $z_{p,q}=z_p-z_q$, $\tilde{z}_{p,q}=\tilde{z}_p-\tilde{z}_q$, $\Delta z_{p}=\tilde{z}_p-z_p$, $\Delta z_{q}=\tilde{z}_q-z_q$ and $\Delta z_{p,q}=\tilde{z}_{p,q}-{z}_{p,q}$. \cref{equ:errbnd} indicates that:
\begin{equation}
\begin{aligned}
\|\Delta z_{p,q}\|_1 & \le \|\Delta z_p\|_1 + \|\Delta z_q\|_1 \\
& \le (\|z_p\|_1depth(p)+\|z_q\|_1depth(q))\epsilon \\
\end{aligned}
\end{equation}  
We can assume that $\|\Delta z_{p,q}\|_1$ is very small and $\Delta z_{p,q}$ is close to $\mathbf{0}$. Then by ignoring the second-order terms, we have:
\begin{equation}
\|\tilde{z}_{p,q}\|_2^2=\|{z}_{p,q}+\Delta z_{p,q}\|_2^2 \approx \|{z}_{p,q}\|_2^2+ 2{z}_{p,q}^T \Delta z_{p,q}~.
\end{equation} 
Let $\tilde{R}_{p,q}=\|\tilde{z}_{p,q}\|_2^2$ denote the approximate effective resistance, then we obtain:
\begin{equation}
\begin{aligned}
|\frac{\tilde{R}_{p,q}}{R_{p,q}}-1| & \approx | \frac{2{z}_{p,q}^T \Delta z_{p,q}}{\|z_{p,q}\|_2^2} | \le \frac{2\|{z}_{p,q}\|_1 \|\Delta z_{p,q}\|_1}{\|z_{p,q}\|_2^2} \\
& \le \frac{2\|{z}_{p,q}\|_1 (\|z_p\|_1depth(p)+\|z_q\|_1depth(q))}{\|z_{p,q}\|_2^2} \epsilon~.\\
\end{aligned}
\end{equation}
If we denote $\alpha_{p,q}=\frac{2\|{z}_{p,q}\|_1 (\|z_p\|_1depth(p)+\|z_q\|_1depth(q))}{\|z_{p,q}\|_2^2}$, then
\begin{equation}
\label{equ:reserrbnd}
1-\alpha_{p,q}\epsilon \le \frac{\tilde{R}_{p,q}}{R_{p,q}} \le 1+\alpha_{p,q}\epsilon~.
\end{equation}
This implies that the relative error of effective resistance scales linearly with the parameter $\epsilon$. The smaller $\epsilon$ is, the closer $\tilde{R}_{p,q}$ is to $R_{p,q}$.

Computing effective resistances using (\ref{equ:approxeffres}) requires Cholesky factorization on $L_G$. However, for large-scale graphs, especially for graphs that stem from social networks, Cholesky factorization on $L_G$ can be very expensive. In this work, we propose to use incomplete Cholesky factorization instead. In incomplete Cholesky factorization, some fill-ins with very small absolute values are dropped, which corresponds to set some branches with large resistances to open and does not introduce large errors to effective resistances. We summarize the overall algorithm for computing effective resistances as follows.
\begin{algorithm}[H]
		\caption{Computing Effective Resistances Based on Sparse Approximate Inverse of Cholesky Factor}
		\label{alg:effres}
		\begin{algorithmic}[1]
			\setcounter{ALC@unique}{0}
			\REQUIRE A weighted undirected graph: $G$, a set of effective resistance queries $Q_r$. 
			\ENSURE Effective resistances for each query in $Q_r$.
			\STATE Run incomplete Cholesky factorization on $L_G$ to obtain: $L_G \approx LL^T$.
			\STATE Compute the sparse and approximate inverse with \cref{alg:iinv} for $L$: $\tilde{Z} \approx L^{-1}$.
			\FOR{each query $(p,q)$ in $Q_r$}
			\STATE Compute effective resistance: $R_G(p,q) \approx \|\tilde{z}_p-\tilde{z}_q\|_2^2$.
			\ENDFOR
		\end{algorithmic}
	\end{algorithm}

The time complexity of the proposed algorithm is closely related to the number of nonzeros in $\tilde{Z}$. In our experiments, the number of nonzeros in $\tilde{Z}$ is about $Cnlogn$ where $C$ is a small constant (e.g. $C<20$). The reason behind this phenomenon may be the decay property of $L^{-1}$ \cite{decay84}. Here we just assume the average number of nonzeros in one column of $\tilde{Z}$ is $O(logn)$. Incomplete Cholesky factorization takes $O(n)$ time and the number of nonzeros in $L$ is $O(n)$. Each iteration of \cref{alg:iinv} takes $O(\frac{nnz(L)}{n}logn+logn \cdot  loglogn)=O(logn \cdot loglogn)$ time and the time complexity of \cref{alg:iinv} is $O(n logn \cdot loglogn)$. So the overall time complexity of \cref{alg:effres} is $O(nlogn \cdot loglogn)+O(|Q_r|logn)$.

\section{Numerical Results}
We first compare the proposed algorithm for computing effective resistances (\cref{alg:effres}) with the algorithm proposed in \cite{SCWWW15}. Because the codes shared by the authors of \cite{SCWWW15} are written in MATLAB, we also implement our \cref{alg:effres} in MATLAB. Then we implement the graph sparsification based power grid reduction (\cref{alg:ResMOR}) and compare the scenarios where effective resistances are computed accurately, approximately using the random projection based method \cite{SCWWW15} and using the proposed \cref{alg:effres}. Finally, the resulted fast power grid reduction algorithm is leveraged to solve problems of DC incremental analysis and transient analysis. These programs are written in C++. All experiments are conducted using a single CPU core of a computer with Intel Xeon E5-2630 CPU @2.40 GHz and 256 GB RAM.

\subsection{Results on Computing Effective Resistances}

In this subsection, we compare the proposed algorithm (\cref{alg:effres}) with the random projection based method \cite{SCWWW15}, whose results are obtained by running the codes shared on \cite{rjeffresweb}. We do not compare the algorithms in \cite{SCIJCAI16,SCKDD21} because these algorithms can only handle unweighted graphs. The results are listed in \cref{table1}. The test cases cover a great variety of graphs obtained from social networks, finite element analysis and circuit simulation problems \cite{thupg,ibmpg,sparse}. For each case, effective resistances of all edges are computed, i.e. the set of effective resistance queries $Q_r=E$.
\begin{table}[h]
		\centering
 		\setlength{\abovecaptionskip}{0pt}
  \caption{Results for Computing Effective Resistances on Large Graphs.}
		\label{table1}
		\resizebox{\linewidth}{!}{
			\begin{tabular}{@{}c@{~}c@{~}c@{~}c@{~}c@{~}c@{~}c@{~}c@{~}c@{~}c@{~}c@{~}c@{}}
				\hline
				\multirow{2}{*}{Case}&\multirow{2}{*}{$|V|(|E|)$} & \multirow{2}{*}{$dpt$} & \multicolumn{4}{c}{WWW15 \cite{SCWWW15}} & & \multicolumn{4}{c}{\cref{alg:effres}} \\
				\cline{4-7} \cline{9-12}
				& & & $T(s)$ & $E_a$ & $E_m$ & $\frac{\mathit{nnz}(Q)}{nlogn}$ & & $T(s)$ & $E_a$ & $E_m$ & $\frac{\mathit{nnz}(\tilde{Z})}{nlogn}$  \\
				\hline
				com-DBLP & 3.2E5(1.0E6)& 464 & 517 & 2.6E-2 & 1.4E-1 & 108 & & 4.14 & 7.1E-5 & 1.9E-3 & 5.40 \\
				\hline
				com-Amaz & 3.3E5(9.3E5) & 590 & 719 & 2.2E-2 & 1.4E-1 & 149 & & 4.71 & 8.0E-5 & 3.9E-3 & 7.47 \\
				\hline
				com-Yout & 1.1E6(3.0E6) & 1370 & 926 & 3.5E-2 & 2.1E-1 & 32.6 & & 21.0 & 1.5E-4 & 2.1E-2 & 1.63  \\
				\hline
				coAuDBLP & 3.0E5(1.0E6) & 1040 & 513 & 2.5E-2 & 1.1E-1 & 108 & & 3.87 & 7.1E-5 & 4.0E-3 & 5.43 \\
				\hline
				coAuCite & 2.3E5(8.1E5) & 774 & 414 & 2.4E-2 & 1.0E-1 & 129 & & 2.32 & 5.6E-5 & 7.9E-3 & 6.45\\
				\hline
				fe\_tooth & 7.8E4(4.5E5) & 1892 & 322 & 1.8E-2 & 7.4E-2 & 304 & & 1.73 & 8.6E-4 & 1.1E-2 & 15.2\\
				\hline
				fe\_rotor & 1.0E5(7.6E5) & 2448 & 488 & 1.7E-2 & 7.0E-2 & 344 & & 2.84 & 8.3E-4 & 2.1E-2 & 17.2 \\
				\hline
				NACA0015 & 1.0E6(3.1E6) & 543 & 2447 & 2.2E-2 & 7.5E-2 & 163 & & 12.1 & 1.0E-3 & 3.6E-3 & 8.17\\
				\hline
				ibmpg5 & 1.1E6(1.6E6) & 513 & 691 & 2.2E-2 & 1.2E-1 & 123 & & 3.16 & 1.7E-3 & 2.7E-2 & 6.17 \\
				\hline
				ibmpg6 & 1.7E6(2.5E6) & 602 & 934 & 2.3E-2 & 1.2E-1 & 109 & & 4.64 & 1.8E-3 & 2.2E-2 & 5.44 \\
				\hline
				thupg1 & 5.0E6(8.2E6) & 1097 & 7158 & 1.8E-2 & 8.1E-2 & 122 & & 36.6 & 1.7E-3 & 1.4E-2 & 6.10\\
				\hline
				G2\_circuit & 1.5E5(2.9E5) & 720 & 214 & 2.0E-2 & 1.2E-1 & 166 & & 1.15 & 1.3E-3 & 4.4E-2 & 8.30 \\
				\hline
				G3\_circuit & 1.6E6(3.1E6) & 1237 & 2388 & 2.0E-2 & 9.8E-2 & 140 & & 13.3 & 3.1E-3 & 4.0E-2 & 7.02\\
				\hline
				thupg10 & 6.0E7(1.0E8) & 3725 & -& - & - & - & & 481 & 1.7E-3 & 1.7E-2 & 1.24\\
				\hline
			\end{tabular}
		}
	\end{table}
$T$ denotes the runtime for computing effective resistances. $E_a$ and $E_m$ denote the average and the maximum relative errors, which are estimated by randomly selecting 1000 edges, computing accurate effective resistances for these edges and then computing the errors. For \cref{alg:effres}, $T$ includes the time for incomplete Cholesky factorization, computing approximate inverse and computing effective resistances. The drop tolerance in incomplete Cholesky factorization is set to 1E-3 and the parameter $\epsilon$ in \cref{alg:iinv} is also set to 1E-3. We also record the maximum depth of the filled graph, which is defined in the last section and denoted by $dpt$, and the number of nonzeros in the random projection matrix $Q$ and in the appproximate inverse matrix $\tilde{Z}$, both of which are divided by $nlogn$. ``-" means that it takes more than 10 hours. 

From the results we see that, the proposed \cref{alg:effres} achieves an average \textbf{168X} speedup over the random projection based method \cite{SCWWW15}. Although the number of nonzeros in the random projection matrix $Q$ is much larger than that in the approximate inverse matrix $\tilde{Z}$, the proposed algorithm shows one to two orders of magnitude improvement in the avarage relative error and also significant reduction in the maximum relative error. For the largest case named ``thupg10", with 6.0E7 nodes and 1.0E8 edges, effective resistances for up to 1.0E8 node pairs can be computed within just 481 seconds, while the average relative error is about 0.17\%, which demonstrates the extremely high scalability of the proposed algorithm.

\subsection{Results on Graph Sparsification Based PG Reduction}

Combining the proposed algorithm for computing effective resistances (\cref{alg:effres}) with the graph sparsification based power grid reduction framework (\cref{alg:ResMOR}), we obtain a fast power grid reduction algorithm. Test cases are from the well-known IBM power grid benchmarks \cite{ibmpg}. For graph partitioning, we use the widely adopted graph partitioner METIS \cite{metis} and the number of blocks are set to $\frac{ \# ports}{50}$. For transient analysis, each case is simulated for $1000$ fixed-size time steps and both original models and reduced models are analyzed with the direct solver CHOLMOD \cite{cholmod2008} (performing just once matrix factorization).  For DC incremental analysis, $10\%$ blocks of each benchmark are modified to mimick the design process where the initial power grid is modified to fix violations. We compare three power grid reduction methods, which only differ in computing effective resistances. The first method computes effective resistances accurately, while the other two compute effective resistances approximately, using the random projection based approach \cite{SCWWW15} and the proposed approach (\cref{alg:effres}) respectively. 

The results are listed in \cref{table2}. $|V|$ and $|E|$ denote the number of nodes and resistors. $T_{\mathit{red}}$ denotes the time for power grid reduction. Note that in the application of DC incremental analysis, only the modified blocks need to be reduced, so the $T_{\mathit{red}}$ in incremental analysis is just about 10\% of that in transient analysis. $T_{\mathit{tr}}$ and $T_{\mathit{inc}}$ are the time for transient analysis and DC incremental analysis, respectively. Err is the average absolute error and Rel is the relative error, which is computed by dividing the Err by the maximum voltage drop.
\begin{table*}[t]
		\centering
 		\setlength{\abovecaptionskip}{0pt}
		\caption{Results on Graph Sparsification Based Power Grid Reduction for Transient Analysis (upper) and DC Incremental Analysis (lower).}
		\label{table2}
		\resizebox{\linewidth}{!}{
			\begin{tabular}{@{~}c@{~}c@{~}c@{~}c@{~}c@{~}c@{~}c@{~}c@{~}c@{~}c@{~}c@{~}c@{~}c@{~}c@{~}c@{~}c@{~}c@{~}c@{~}c@{~}c@{~}c@{~}}
				\hline
				\multirow{2}{*}{Case}& \multicolumn{2}{c}{Original} & & \multicolumn{5}{c}{w/ Acc. Eff. Res.} & & \multicolumn{5}{c}{w/ App. Eff. Res. (\cite{SCWWW15})} & & \multicolumn{5}{c}{w/ App. Eff. Res. (\cref{alg:effres})} \\
				\cline{2-3} \cline{5-9} \cline{11-15} \cline{17-21}
				& $|V|(|E|)$ & $T_{\mathit{tr}}$ & & $|V|(|E|)$ & $T_{\mathit{red}}$ & $T_{\mathit{tr}}$ & Err(mV) & Rel(\%) & & $|V|(|E|)$ & $T_{\mathit{red}}$ & $T_{\mathit{tr}}$ & Err(mV) & Rel(\%) & & $|V|(|E|)$ & $T_{\mathit{red}}$ & $T_{\mathit{tr}}$ & Err(mV) & Rel(\%)   \\
				\hline
				ibmpg2t & 1.3E5(2.08E5) & 13.3 & & 5.1E4(1.49E5) & 6.55 & 4.33 & 0.801 & 1.52 & & 5.1E4(1.47E5) & 3.51 & 4.23 & 2.300 & 4.28 & & 5.1E4(1.49E5) & 0.951 & 4.20 & 0.796 & 1.51 \\
				\hline
				ibmpg3t & 8.5E5(1.40E6) & 201 & & 3.0E5(8.59E5) & 67.2 & 35.1 & 0.153 & 0.78 & & 3.0E4(8.54E5) & 29.4 & 34.7 & 0.253 & 1.29 & & 3.0E4(8.59E5) & 7.70 & 35.8 & 0.162 & 0.83 \\
				\hline
				ibmpg4t & 9.5E5(1.55E6) & 310 & & 4.0E5(1.35E6) & 81.9 & 132 & 0.006 & 0.93 & & 4.0E5(1.35E6) & 36.3 & 133 & 0.034 & 4.85 & & 4.0E5(1.35E6) & 10.6 & 129 & 0.006 & 0.93 \\
				\hline
				ibmpg5t & 1.1E6(1.62E6) & 141 & & 5.1E5(1.04E6) & 24.1 & 48.0 & 0.124 & 0.87 & & 5.1E5(1.03E6) & 21.5 & 47.4 & 0.137 & 0.96 & & 5.1E5(1.04E6) & 5.59 & 48.8 & 0.125 & 0.87 \\
				\hline
				ibmpg6t & 1.7E6(2.48E6) & 179 & & 7.7E5(1.64E6) & 39.4 & 67.7 & 0.173 & 1.02 & & 7.7E5(1.66E5) & 33.5 & 67.0 & 0.334 & 1.97 & & 7.7E5(1.64E6) & 8.76 & 67.5 & 0.173 & 1.02 \\
				\hline				
			\end{tabular}
		}
\resizebox{\linewidth}{!}{
			\begin{tabular}{@{~}c@{~}c@{~}c@{~}c@{~}c@{~}c@{~}c@{~}c@{~}c@{~}c@{~}c@{~}c@{~}c@{~}c@{~}c@{~}c@{~}c@{~}c@{~}c@{~}c@{~}c@{~}}
				\hline
				\multirow{2}{*}{Case}& \multicolumn{2}{c}{Original} & & \multicolumn{5}{c}{w/ Acc. Eff. Res.} & & \multicolumn{5}{c}{w/ App. Eff. Res. (\cite{SCWWW15})} & & \multicolumn{5}{c}{w/ App. Eff. Res. (\cref{alg:effres})} \\
				\cline{2-3} \cline{5-9} \cline{11-15} \cline{17-21}
				& $|V|(|E|)$ & $T_{\mathit{inc}}$ & & $|V|(|E|)$ & $T_{\mathit{red}}$ & $T_{\mathit{inc}}$ & Err(mV) & Rel(\%) & & $|V|(|E|)$ & $T_{\mathit{red}}$ & $T_{\mathit{inc}}$ & Err(mV) & Rel(\%) & & $|V|(|E|)$ & $T_{\mathit{red}}$ & $T_{\mathit{inc}}$ & Err(mV) & Rel(\%)  
				\\
				\hline
				ibmpg2 & 1.3E5(2.08E5) & 0.627 & & 5.1E4(1.49E5) & 0.784 & 0.159 & 6.052 & 1.21 & & 5.1E4(1.47E5) & 0.381 & 0.149 & 14.10 & 2.81 & & 5.1E4(1.49E5) & 0.115 & 0.158 & 6.020 & 1.20 \\
				\hline
				ibmpg3 & 8.5E5(1.40E6) & 20.4 & & 3.0E5(8.59E5) & 6.73 & 1.54 & 1.612 & 0.66 & & 3.0E4(8.54E5) & 2.66 & 1.46 & 3.004 & 1.22 & & 3.0E4(8.59E5) & 0.769 & 1.57 & 1.734 & 0.71 \\
				\hline
				ibmpg4 & 9.5E5(1.55E6) & 36.3 & & 4.0E5(1.35E6) & 8.85 & 11.7 & 0.093 & 0.76 & & 4.0E5(1.35E6) & 3.94 & 11.5 & 1.159 & 9.45 & & 4.0E5(1.35E6) & 1.04 & 11.6 & 0.089 & 0.73 \\
				\hline
				ibmpg5 & 1.1E6(1.62E6) & 10.5 & & 5.1E5(1.04E6) & 2.32 & 1.78 & 0.929 & 1.34 & & 5.1E5(1.03E6) & 2.11 & 1.74 & 5.257 & 7.59 & & 5.1E5(1.04E6) & 0.538 & 1.81 & 0.929 & 1.34 \\
				\hline
				ibmpg6 & 1.7E6(2.48E6) & 8.68 & & 7.7E5(1.64E6) & 3.98 & 2.42 & 1.503 & 0.73 & & 7.7E5(1.66E5) & 3.53 & 2.45 & 3.181 & 1.54 & & 7.7E5(1.64E6) & 0.888 & 2.41 & 1.512 & 0.73 \\
				\hline				
			\end{tabular}
		}
	\end{table*}

From the results we see that, with the proposed algorithm for computing effective resistances (\cref{alg:effres}), the time for power grid reduction can be reduced largely, showing an average 6.4X speedup over the power grid reduction approach based on accurate effective resistances.
In terms of the reduction accuracy, there is almost no increase in reduction errors. This validates the efficiency and the effectiveness of the proposed power grid reduction approach. On the other hand, with the random projection based method, it still takes a large amount of time to generate reduced models. Meanwhile, the reduction accuracy is also deteriorated, which is due to the large errors of effective resistances.

\begin{figure}[h]
		\setlength{\abovecaptionskip}{0pt}
		\setlength{\belowcaptionskip}{0pt}
		\centering
		\includegraphics[width=9.0cm]{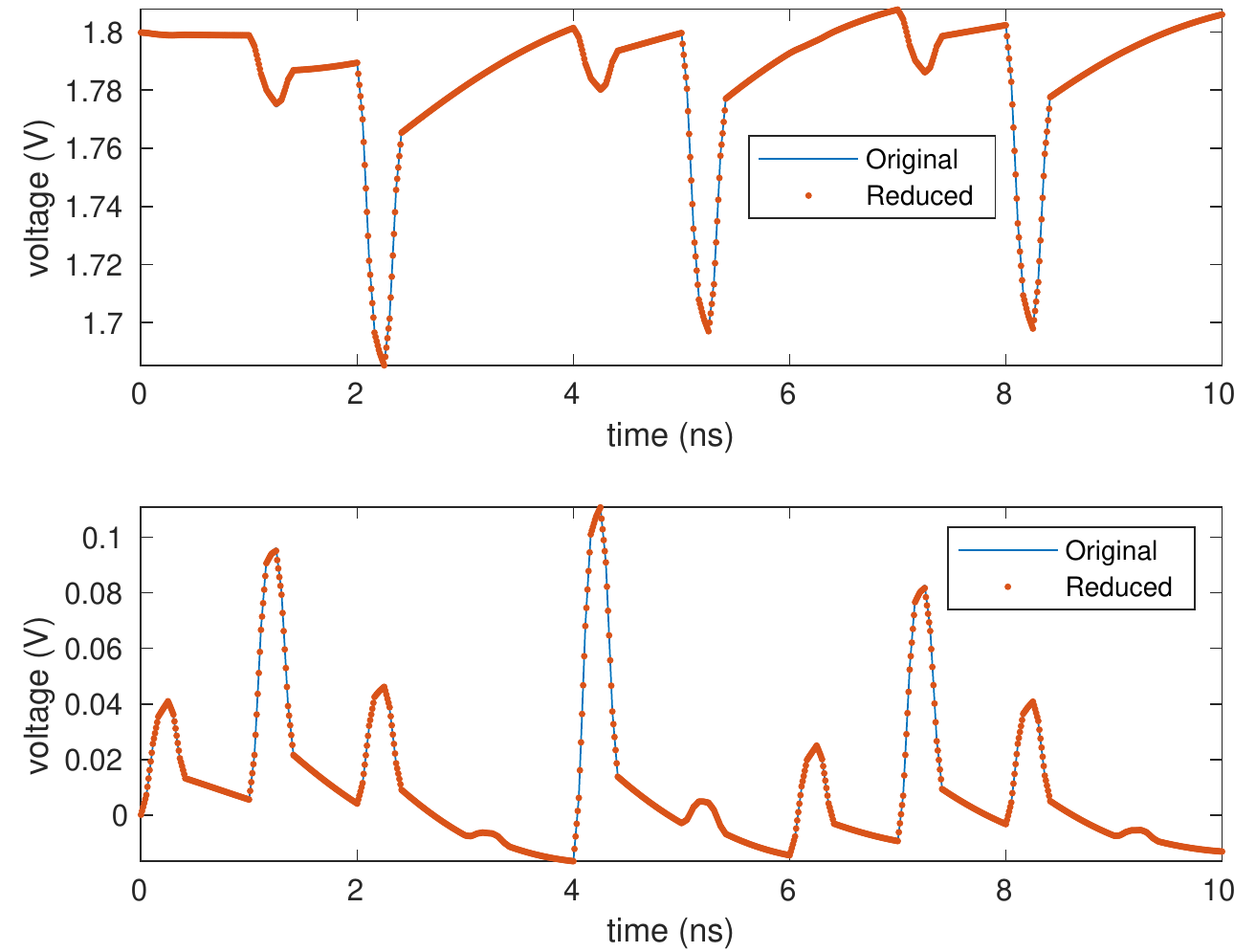}
		\caption{The transient simulation results of a VDD node (up) and a GND node (down) in case ``ibmpg3t'', obtained by analyzing the original power grid and the reduced power grid.}
		\label{fig:waveform}
	\end{figure}

As for the total time of transient analysis, which includes the time for power grid reduction and the time for analyzing the reduced power grid, the proposed method achieves 1.7X speedup averagely over the power grid reduction approach based on accurate effective resistances. Compared with analyzing the original power grid, the proposed method shows an average 2.9X speedup, with the relative error below 1.51\% for all cases. The transient waveforms of node n0\_20706300\_8937900 and node n1\_29561400\_9521100 in case ``ibmpg3t" are plotted in \cref{fig:waveform}. They validate the accuracy of transient simulation using the proposed fast power grid reduction method. 

In the application of DC incremental analysis, the proposed method shows significant improvement in the incremental power grid reduction stage, thus achieves an average 2.5X speedup for the total time over the power grid reduction method based on accurate effective resistances. Compared with analyzing the modified power grid directly using CHOLMOD ("Original" in Table II), the proposed method shows 4.2X speedups on average, while the relative error is below 1.34\% for all cases. 

\section{Conclusions}
In this paper, we propose an efficient algorithm for computing effective resistances on massive weighted graphs. The proposed algorithm is based on a sparse approximate inverse technique for Cholesky factors and achieves 168X speedups on average and also significant reduction in errors, compared with the random projection based method. Combining the proposed algorithm with the graph sparsification based power grid reduction framework, we obtain a fast power grid reduction method, which shows 6.4X speedups averagely and almost the same reduction accuracy. Utilized for the downstream applications of power grid transient analysis and incremental analysis, the proposed fast power grid reduction method brings 1.7X and 2.5X speedups of overall time, compared to using the reduction method based on accurate effective resistances. Given the wide range of applications of effective resistances, we hope that this work can promote progress in more fields in the future.   

\bibliography{ref1}
\bibliographystyle{IEEEtran}

\end{document}